\documentclass[10pt]{amsart}
\usepackage{caption}
\usepackage{subcaption}
\usepackage{enumerate,url,amssymb,  mathrsfs, graphicx, pdfsync}
\newtheorem{theorem}{Theorem}[section]

\newtheorem*{lemma*}{Lemma}

\newtheorem{proposition}[theorem]{Proposition}

\theoremstyle{definition}

\newtheorem{question}[theorem]{Question}

\theoremstyle{remark}

\numberwithin{equation}{section}

\newcommand{\abs}[1]{\lvert#1\rvert}

\newcommand{\norm}[1]{\lVert#1\rVert}

\newcommand{\A}{\mathbb{A}}
\newcommand{\C}{\mathbb{C}}

\newcommand{\W}{\mathscr{W}}

\newcommand{\R}{\mathbb{R}}

\newcommand{\dtext}{\textnormal d}

\newcommand{\onto}{\xrightarrow[]{{}_{\!\!\textnormal{onto\,\,}\!\!}}}

\newcommand{\Hp}[1]{\mathcal H^{1,#1} (\A, \A^\ast)}
\newcommand{\Hpcl}[1]{\overline{\mathcal H}^{1,#1} (\A, \A^\ast)}
\newcommand{\Hid}[1]{\mathcal H_{\Id}^{1,#1} (\A, \A^\ast)}

\newcommand{\Rp}[1]{\mathcal R^{1,#1} (\A, \A^\ast)}
\newcommand{\tRp}[1]{\overline{\mathcal R}^{1,#1} (\A, \A^\ast)}

\newcommand{\ee}{\mathbb{E}}

\DeclareMathOperator{\im}{Im}

\DeclareMathOperator{\loc}{loc}

\DeclareMathOperator{\Diff}{Diff}
\DeclareMathOperator{\Mod}{Mod}

\DeclareMathOperator{\Id}{Id}
\def\le{\leqslant}
\def\ge{\geqslant}

\def\XXint#1#2#3{{\setbox0=\hbox{$#1{#2#3}{\int}$}\vcenter{\hbox{$#2#3$}}\kern-.5\wd0}}

\def\XXiint#1#2#3{{\setbox0=\hbox{$#1{#2#3}{\iint}$}\vcenter{\hbox{$#2#3$}}\kern-.5\wd0}}

\begin{document}
\title{Radially symmetry of minimizers to the weighted Dirichlet energy}

\author[A. Koski]{Aleksis Koski}
\address{Department of Mathematics and Statistics, P.O.Box 35 (MaD) FI-40014 University of Jyv\"askyl\"a, Finland}
\email{aleksis.t.koski@jyu.fi}

\author[J. Onninen]{Jani Onninen}
\address{Department of Mathematics, Syracuse University, Syracuse,
NY 13244, USA and  Department of Mathematics and Statistics, P.O.Box 35 (MaD) FI-40014 University of Jyv\"askyl\"a, Finland
}
\email{jkonnine@syr.edu}
\thanks{  A. Koski was supported by the Academy of Finland Grant number 307023.
J. Onninen was supported by the NSF grant  DMS-1700274.}





\maketitle

\begin{abstract} We consider the problem of minimizing the weighted Dirichlet energy between homeomorphisms of planar annuli.  A known challenge lies in the  case  when the weight $\lambda$ depends on the independent variable $z$.  We prove that for an increasing radial weight $\lambda(z)$  the infimal energy within the class of all Sobolev homeomorphisms is the same as in the class of radially symmetric maps. For a general radial  weight $\lambda(z)$  we establish the same result in the case when the target is conformally thin compared to the domain. Fixing the admissible homeomorphisms  on the outer boundary we establish  the radial symmetry for every such weight.

\end{abstract}
\section{Introduction}

The Geometric function theory (GFT) is largely concerned with generalizing the theory of conformal mappings. In this paper we turn
to a variational approach and study  mappings with
smallest weighted mean distortion. Our underlying idea is to  extend the theory of extremal quasiconformal mappings to the  minima of the weighted  mean average   of the distortion function; that is, minimize  the weighted $\mathscr L^1$-norm as opposed to the $\mathscr L^\infty$-norm in the Teichm\"uller theory.    There are many natural reasons for studying  such a minimization problem.  This quickly leads one to    extremal mappings of mean distortion  between annuli, a classical and well-understood problem for extremal quasiconformal mappings, traditionally referred to  as the Gr\"otsch problem.  Indeed, annuli
\[\A=A(r,R)= \{z\in \R^n \colon r< \abs{z} <R\}  \quad \textnormal{ and } \quad \A^\ast=\{z\in \R^n \colon r_\ast< \abs{z} <R_\ast\}  \]
 are the first in order of complexity where one observes nontrivial conformal invariants such as moduli. Hereafter $0< r< R$ and $0<  r_\ast <R_\ast$ are called the inner and outer radii of $\A$ and $\A^\ast$, respectively.  Precisely, we search mappings with  least weighted distortion 
 \begin{equation}\label{eq:Kmini}
  \int_{\A^\ast} \lambda (|z|, |f|) K_f(z) \, \dtext z
  \end{equation}
 subject to Sobolev homeomorphims $f \colon \A^\ast  \onto \A$ with integrable distortion.
 Recall that a Sobolev homeomorphism  $f \in \W_{\loc}^{1,1} (\A^\ast, \C)$ has  finite distortion if  there is a measurable function $K (z) \ge 1$, finite a.e., such that
\begin{equation}\label{eq:fd}
\abs{Df(z)}^2 \le  2\, K(z) \, J_f(z)\, , \qquad J_f(z)=\det Df(z) .
\end{equation}
Here and in what follows we use  the Hilbert-Schmidt norm of  a linear map $A$, defined by the rule $ \norm{A}^2 =  \langle A\,,\,A \rangle \,=\,\textnormal{Tr} [A^* \cdot A ]$. We are interested in homeomorphisms and their limits, we recall that the Jacobian determinant of a Sobolev homeomorphism is always locally integrable.  The distortion inequality~\eqref{eq:fd} merely asks that  the differential $Df(z)$ vanishes at those points $z$ where the Jacobian $J_f(z)=0$. . 
The smallest function $K(z) \geqslant 1$ for which the distortion inequality~\eqref{eq:fd} holds is denoted by $K_f(z)$,
\[
K_f (z) = \begin{cases} \frac{\abs{Df(z)}^2 }{2 \, J_f(z)} \quad & \textnormal{if } J (z,f) >0 \\
1 & \textnormal{if } J (z,f) =0  \end{cases}
\]
We obtain quasiconformal mappings $f$ if $K\in \mathscr L^\infty (\A^\ast)$.  The theory of quasiconformal mappings is by now well understood, see the monographs \cite{Reb} by Reshetnyak, \cite{Rib} by Rickman and \cite{IMb} by Iwaniec and Martin.  In the last 20 years, systematic studies of mappings of finite distortion have emerged in GFT~\cite{AIMb, IMb, HKb}.  The theory of mappings of finite distortion arose from the  need to extend the ideas and applications of the classical theory of quasiconformal mappings to the degenerate elliptic setting.  Motivated by mathematical models of nonlinear elasticity~\cite{Ba, Ba12, Ba2}, the focus has been finding a class of mappings, as close to homeomorphisms as possible, in which the minimum energy~\eqref{eq:Kmini} is attained. In the case of minimizing the weighted  $\mathscr L^1$-mean distortion this is  possible only when we move, equivalently, to minimize the weighted  Dirichlet energy of the
inverse map. Indeed,  the inverse map  $h = f^{-1} : \A \onto \A^\ast$ belongs to the Sobolev space $\mathscr W^{1,2}(\A ,\mathbb R^2)$, and we have
 \begin{equation}\label{eq:trans}
 \int_\A \lambda (\abs{h(z)}, \abs{z} )|Dh(z)|^2 \,\dtext z  = 2\, \int_{\A^\ast }  \lambda (\abs{z}, \abs{f(z)} ) K_f(z)\,\dtext  z \, , 
 \end{equation}
 see~\cite{AIMO, HK, HKO} for details. 
 In such a problem  the transition to the energy of the inverse map  results in a convex variational integral. Therefore, from now on we  minimize the weighted Dirichlet energy
\begin{equation}\label{eq:minipr}
 \mathbb E_\lambda [h]=  \int_\A \lambda (\abs{h(z)}, \abs{z} )  \abs{Dh(z)}^2 \, \dtext z \, .
 \end{equation}
The infimum is subjected to orientation preserving Sobolev homeomorphisms $h \colon \A $ $\onto \A^\ast$ in $\W^{1,2} (\A, \C)$ which are furthermore assumed to preserve the order of the boundary components. Such a class of Sobolev homeomorphisms is denoted by $\Hp{2}$.  

Because of rotational symmetry it seems likely that the energy-minimal deformations of~\eqref{eq:minipr} are radial minimal mappings. However, the difficulty to verify the rotationally symmetry is well recognized in the theory of nonlinear elasticity.   A number of papers in the literature is devoted to understand the expected radial symmetric properties~\cite{AIM, Ba1, CK, CG, HLW, Ho,  IKOni, IKO3,  IOhy, IOne, IOan, JM, JK,  Ka, Ka2, KOpharm, Me, MS, S, SiSp, St}. We study this question for the weighted Dirichlet energy.  
\begin{question}\label{qu:main}
Does the equality
\begin{equation}\label{eq:mainq}
\underset{\Hp{2}}{\inf} \mathbb E_\lambda [h] = \underset{\mathcal R^{1,2} (\A, \A^\ast)}{\inf}   \mathbb E_\lambda [h]  
\end{equation}
hold?
\end{question}
 In what follows, we denote the subclass of radial homeomorphisms by
\[\mathcal R^{1,2} (\A, \A^\ast) = \left\{ h \in \Hp{2} \colon h(x)= H(\abs{z}) \frac{z}{\abs{z}}  \right\} \, .\] 
In general the mappings with the infimum energy in~\eqref{eq:mainq} need not be homeomorphisms. In fact a part of the domain near its boundary may collapse into the inner boundary of the target annulus $\A^\ast$. In mathematical models of nonlinear elasticity this is interpreted as {\it interpenetration of matter}. Of course, in general enlarging the set of the admissible
mappings may change the nature of the energy-minimal solutions. This may result in smaller energy than in $\Hp{2}$, whether or not the infimum is attained. To avoid such an effect one needs to know that a weak limit of a minimizing sequence $h_j \in  \Hp{2}$ can be realized as a strong limit of homeomorphisms in $\W^{1,2} (\A, \C)$. This follows from the result in~\cite{IOw=s}  which tells us that the classes of the weak
and strong limits of $\W^{1,2}$-Sobolev homeomorphisms are the same. We denote such a class of deformations from $\A$ into $\overline{\A^\ast}$ by $\Hpcl{2}$. It is quite easy to see that mappings in $\Hpcl{2}$ extend as continuous monotone maps of $\overline {\A }$ onto $\overline{\A^\ast}$. As a converse Iwaniec and Onninen~\cite{IOmo} proved a Sobolev
variant of the classical Youngs approximation theorem. According to their result the class $\Hpcl{2}$ equals the class of  orientation-preserving monotone mappings from $\overline{\A}$ onto  $\overline{\A^\ast}$ in the Sobolev class $\W^{1,2} (\A, \C)$ that also preserve the order of the boundary components of annuli. Monotonicity, the concept of C.B. Morrey~\cite{Mor} simply means that the preimage of a point in $ \overline{\A^\ast}$ is a continuum in $\overline{\A}$.

The class of the weak $\W^{1,2}$-limits of  radially symmetric homeomorphisms can be also give purely analytical  characterization.   This leads us to define  
\[  \tRp{2} = \{  h \in \W^{1,2} (\A, \R^n) \colon h(x)= H(\abs{x}) \frac{x}{\abs{x}} \textnormal{ and }  H \in \mathscr R \}\]
where
\[\mathscr R = \{ H \colon [r,R] \to [r_\ast , R_\ast] \colon  H(r)=r_\ast, \,   H(R)=R_\ast  \textnormal{ and } \dot{H}\ge 0 \} \, . \]

Now, the identity~\eqref{eq:mainq} in Question~\ref{qu:main} can be equivalently written as
\[ \underset{\Hpcl{2}}{\min} \mathbb E_\lambda [h]  =  \underset{\tRp{2}}{\min}   \mathbb E_\lambda [h] \, . \]

 First, let us look at this question in  the simplest Dirichlet energy case, $\lambda \equiv 1$. The papers~\cite{AIM, IOan} introduced  and heavily relied the concept of {\it Free Lagrangians} and gave a positive answer to Question~\ref{qu:main} in terms of the Dirichlet energy.  The underlying idea was to estimate the
integrand of an energy functional from below  in terms of free Lagrangians. Free Lagrangians, special null Lagrangians~\cite{Bac},  are nonlinear differential forms whose integral means depends only on  a given homotopy class of homeomorphisms.  The volume form is not only a simple  example of free Lagrangian but also a key player in the proof. In the conformally equivalent  case,  it is the only free Lagrangian that is needed. In spite of being a trivial case, it apprehends the essence
of free Lagrangians well. Suppose that $\A=\A^\ast$ and $h \in \Hp{2}$.   Showing that the identity map  $\Id (z) =z$ is a global minimizer follows immediately from Hadamard's inequality $\abs{Dh (z)}^2 \ge 2  J_h (z)$:
\begin{equation}\label{eq:confgen}
\mathbb E_1[h]  = \int_\A \abs{Dh(z)}^2 \,  \dtext z \ge  2\, \int_\A J_h (z) \, \dtext z  = 2 \int_{\A^\ast} 1\,  \dtext y =  \mathbb E_1[\Id] \, . 
\end{equation} 
Second,  having the techniques of free Lagrangians and the estimates for the Dirichlet energy in hand  proving corresponding results for  weights $\lambda$ that depend only on $|h|$ is  straightforward. Again, the conformal case, captures  this phenomenon well.  Indeed, for $\lambda (\abs{h}, \abs{z}) = \lambda( {\abs{h}}) $, $\A=\A^\ast$ and $h \in \Hp{2}$ we have
 \[
 \mathbb E_\lambda [h]  = \int_\A  \lambda (\abs{h})\, \abs{Dh(z)}^2 \,  \dtext z \ge  2\, \int_\A  \lambda (\abs{h})\, J_h (z) \, \dtext z = 2 \int_{\A^\ast} \lambda (\abs{y}) \, \dtext y  = \mathbb E_\lambda[\Id] \, . 
 \]  
 It is known that the radial minimizers to the weighted Dirichlet energy are the absolute minimizers provided the weight  is independent of  $|z|$, ~\cite{AIM, IOan, Ka} 
\begin{proposition}\label{pr:dirichlet}
Suppose that  $\A , \A^\ast \subset \C$ and $\lambda (\abs{h}, \abs{z}) = \lambda ({\abs{h}})$. Then  the equality~\eqref{eq:mainq} always  holds.
\end{proposition} 
We hence turn our attention in Question~\ref{qu:main} to the case where the weight depends on $\abs{z}$. From now on we assume that the weight $\lambda$ has  the form $\lambda = \lambda (\abs{z})$.  Such a weight brings a completely new challenge to the  studied question. First, there is no trivial case being analogous  to~\eqref{eq:confgen}. Second, such difficulty is already recognized in the literature. The paper~\cite{IOR} is devoted to study the radial minimizers  of $\mathbb E_\lambda $,  $\lambda = \lambda (\abs{z})$.  Question~\ref{qu:main} for the weight  $\lambda = \lambda (\abs{z})$ is explicitly raised in~\cite[Question 4.1]{IOR}.  Our next result proves that the radial minimizers are indeed absolute minimizers when the weight is increasing.

 \begin{theorem}\label{increasingThm}
 Assume that $\lambda \colon [r,R] \to \R$ is continuous, positive and nondecreasing. Then the weighted Dirichlet energy
\begin{equation}\label{weightedDiri}
\ee_{\lambda}[h] = \int_{\A} |Dh(z)|^2 \lambda(\abs{z}) \, \dtext z
\end{equation}
admits a radially symmetric minimizer in the class $\Hpcl{2}$. Moreover, there exists an increasing function $m_\lambda : (0,\infty) \to (0,\infty)$ so that this minimizer is a homeomorphism exactly when $R_*/r_* \geq m_\lambda(R/r)$.
\end{theorem}

For more general weights $\lambda$ the question whether the minimizer of $\ee_{\lambda}$ is attained for a radial mapping remains open in general. However, if the target annulus is conformally thin enough we are able to establish the radial symmetry of the minimizer with no extra assumptions on the weight.
\begin{theorem}\label{thinThm} Let $\lambda \colon [r,R] \to (0, \infty)$ be a continuous function. Then there exists a function $g_\lambda : (1,\infty) \to (1,\infty)$ so that whenever $R_*/r_* \leq g_\lambda(R/r)$, the equality~\eqref{eq:mainq} holds.
\end{theorem}  

 Of course, even in the case of the Dirichlet energy, $\lambda \equiv 1$, the minimizers need not be harmonic. In general, the Euler-Lagrange equation of $\mathbb E_\lambda$ is unavailable; one cannot perform first variations $h+\varepsilon \varphi$ within the class of Sobolev homeomorphisms, not even in $\Rp{2}$.  Therefore, narrowing the admissible homeomorphisms in $\Hp{2}$ does not change the difficulty of the question in this respect.

 \subsection{Partially fixed boundary value problem} Finally, we study the minimization of the weighted Dirichlet energy under mappings fixed on the outer boundary of $\A$, but allowed to be  free on the inner boundary. For simplifying  the notation, we write $\partial_\circ \A = \{x \in \R^n \colon \abs{x}=R\}$ and
\[\Hid{2}  = \left\{h \in \Hp{2}  \colon h \textnormal{ is continuous up to } \partial_\circ \A \textnormal{ and } h(x)= \frac{R_\ast}{R}x \right\} \, .  \]
We prove that keeping the homeomorphisms fixed  in the minimizing sequence on the outer boundary leads the hunted radial symmetry property for an arbitrary weight.
\begin{theorem}\label{fixedThm}
 Assume that $\lambda \colon [r,R] \to (0, \infty)$ is continuous. Then we have
\begin{equation}
\underset{\Hid{2}}{\inf} \int_\A \abs{Dh(z)}^2 \lambda (\abs{z}) \, \dtext z \ =  \underset{\Rp{2}}{\inf} \int_\A \abs{Dh(z)}^2  \lambda (\abs{z}) \, \dtext z \, . 
\end{equation}
  \end{theorem}
The proof of Theorem~\ref{fixedThm} relies on our recent developments for the partially fixed boundary value problem in \cite{KOpharm}.

\section{Analysis of the radial minimizer and the definition of the function $m_\lambda$.}\label{radialsection}

In the paper~\cite{IOR} the radial minimization problem for the energy \eqref{weightedDiri} was completely characterized. It was shown that when subjected to minimization in the class of radial maps $\tRp{2}$, the energy \eqref{weightedDiri} admits a unique minimizer for every continuous positive weight $\lambda \colon [r,R] \to \R$. We will henceforth denote this radial minimizer by $h_0$. Since $h_0$ is radial, we may write $h_0(s e^{i\theta}) = H(s) e^{i\theta}$ for some nondecreasing surjection $H:[r,R] \to [r_*,R_*]$. In~\cite{IOR}, it was shown that the function $H$ may be recovered from $\lambda$ via the following first-order ODE.
\begin{equation}\label{hdiff}
s\lambda(s)\dot{H}(s) = H(s) \Phi(s) \qquad \text{ equivalently } \qquad \frac{\dot{H}}{H} = \frac{\Phi}{s\lambda}.
\end{equation}
Here the function $\Phi: [r,R] \to \R$ is defined by the equation
\begin{equation}\label{phidiff}
\lambda^2(s) - \Phi^2(s) = s \lambda(s) \dot{\Phi}(s),
\end{equation}
at least at the points where the solution $\Phi$ takes nonnegative values. One of the reasons for why we exclude the negative values of $\Phi$ is that otherwise \eqref{hdiff} would imply that $H$ is decreasing at these points. To make a rigorous definition of the function $\Phi$, we first make the following observation, which is an easy consequence from \eqref{phidiff}.
\\\\
\emph{Observation 1.} Any solution $\Phi$ of equation \eqref{phidiff} is increasing at the points $s$ where $\Phi(s) < \lambda(s)$ and decreasing when $\Phi(s) > \lambda(s)$.
\\\\
Since $\lambda$ is positive everywhere, this observation implies that any solution of \eqref{phidiff} is increasing at points where it takes the value zero. Hence every solution has at most one zero, and if such a point $r_0$ exists then the solution is negative on the interval $[r,r_0)$ and positive on $(r_0,R]$. This motivates us to define the function $\Phi$ as follows.

Given the radii $r,R$ and an initial value $\varphi_0$ which may be any real number, we let $\tilde{\Phi}$ be a solution of \eqref{phidiff} on the interval $[r,R]$ with initial data $\tilde{\Phi}(r) = \varphi_0$. The existence and uniqueness of such a solution will follow from the classical ODE theory as soon as we show that the map
\[(s,\Phi) \to \frac{\lambda^2(s) - \Phi^2}{s \lambda(s)}\] is Lipschitz-continuous with respect to the variable $\Phi$. Since $\lambda(s)$ and $s$ are bounded away from zero and infinity, we need only to verify that no solution of \eqref{phidiff} may blow up. But this follows easily from Observation 1, as the observation implies that any solution is bounded by the number $\max(|\varphi_0|,\max_s \lambda(s))$. Note also that by uniqueness the graphs of any two solutions to \eqref{phidiff} do not intersect.

We now define $\Phi$ by $\Phi = \max\{0,\tilde{\Phi}\}$. From the discussion after Observation 1 we know that either $\Phi(s) = \tilde{\Phi}(s)$ everywhere (when $\varphi_0 \geq 0$) or there exists a point $r_0 \in (r,R]$ so that $\Phi \equiv 0$ on $[r,r_0]$ and $\Phi(s) = \tilde{\Phi}(s)$ on $[r_0,R]$ (when $\varphi_0 < 0$).

Given the function $\Phi$, one may always solve the separable ODE \eqref{hdiff} to obtain the function $H$. Furthermore, the conformal modulus of the target annulus $\A^*$ is related to $\Phi$ by the equation
\begin{equation}\label{modulusEquation}
\Mod \A^\ast := \log(R_*/r_*) = \int_r^R \frac{\dot{H}(s)}{H(s)} \, ds = \int_r^R \frac{\Phi(s)}{s \lambda(s)} \, ds.
\end{equation}
Hence the target annulus is defined, up to a scale, by the choice of the initial value $\varphi_0$. As $\varphi_0$ goes through every real number, the equation \eqref{modulusEquation} and the definition of $\Phi$ show that the conformal modulus of the target, $\Mod \A^\ast$, takes every value from $0$ to $\infty$. Hence every possible target annulus is covered by this consideration. We note also that increasing the initial value $\varphi_0$ increases both the function $\Phi$ and the conformal modulus of the target $\A^*$, i.e. makes the target annulus thicker.

\begin{figure}[h]
\centering
\begin{subfigure}{.5\textwidth}
  \centering
  \includegraphics[width=\linewidth]{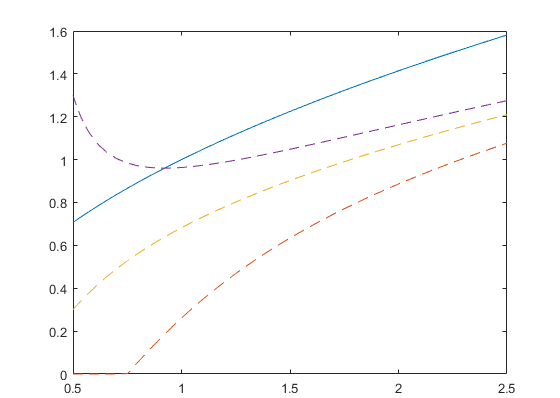}
  \label{fig:sub1}
\end{subfigure}%
\begin{subfigure}{.5\textwidth}
  \centering
  \includegraphics[width=\linewidth]{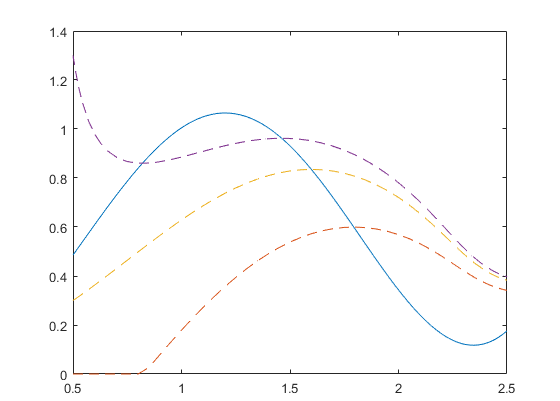}
  \label{fig:sub2}
\end{subfigure}
\caption{Two examples on the annulus $\A(0.5,2.5)$, where the weight $\lambda$ is denoted by a solid line and three different initial values $\varphi_0$ are chosen. The three possible instances of $\Phi$ are shown by the dashed lines, and each of them give rise to a different target annulus.}
\label{fig:test}
\end{figure}
Furthermore, we define the function $m_\lambda$ as follows. Given the numbers $r$ and $R$, we let $\Phi_0$ denote the solution of the equation \eqref{phidiff} with the initial value $\varphi_0 = 0$. We then define $m_\lambda$ by
\begin{equation}\label{mlambda}
\log m_\lambda(R/r) = \int_r^R \frac{\Phi_0(s)}{s \lambda(s)} \, ds.
\end{equation}
Since both equation \eqref{phidiff} and the integral on the right hand side in the above equation are independent with respect to scaling in $s$, the function $m_\lambda$ indeed only depends on the quotient $R/r$ rather than both $r$ and $R$. The above discussion now shows that for every target annulus $\A^*$ with $R_*/r_* \geq m_\lambda(R/r)$ the initial data $\varphi_0$ is nonnegative and hence the map $\Phi$ is positive on $(r,R]$. Likewise for every target $\A^*$ with $R_*/r_* < m_\lambda(R/r)$ there is some radius $r_0 > r$ for which $\Phi \equiv 0$ on $[r,r_0]$. Comparing with \eqref{hdiff}, we find that the radial minimizer $h_0$ is a homeomorphism exactly when $R_*/r_* \geq m_\lambda(R/r)$ (since $\dot{H}(s) > 0$ for $s \in (r,R]$) and for thinner targets $\A^*$ the map $h_0$ fails to be a homeomorphism on the subset $\A(r,r_0) \subset \A(r,R)$ which is collapsed onto the inner boundary of the target (since $\dot{H}(s) = 0$ for $s \in [r,r_0]$).

\section{Free Lagrangians}\label{FLsection}
 In 1977 a novel approach towards minimization of polyconvex energy
functionals was developed and published by J. M. Ball~\cite{Bac}. The
underlying idea was to view the integrand of an energy functional as convex function of null
Lagrangians. The term null Lagrangian pertains to a nonlinear
differential expression whose integral over any open region depends
only on the boundary values of the map,  see~\cite{BCO, Ed, Iw}. Our homeomorphisms $h \colon
\mathbb A \onto {\mathbb A}^\ast $ are not prescribed on the boundary.  
There still exist some nonlinear differential forms, called  {\it free Lagrangians}~\cite{IOan}, defined on a given homotopy class of 
homeomorphisms, whose integral means remain independent of the mapping. These
are rather special null Lagrangians.

Let $\mathbb A = A(r,R)$ and $\mathbb A^\ast = A(r_\ast,R_\ast)$ be two circular annuli in $\C$.   Recall here we  work with one particular homotopy class $\Hp{2}$ of $\W^{1,2}$--orientation-preserving homeomorphisms $ h \colon \A \onto \A^\ast$. that  also preserve the order of the boundary components; that is, $|h(z)| = r_\ast$ for $|z|=r$ and $|h(z)| = R_\ast$ for $|z|=R$.

Clearly, the polar coordinates
\begin{equation}\label{eq:polar}
z= s e^{i\theta} \, , \qquad  r< s <R \quad  \textnormal{ and } \quad  0 \leqslant \theta < 2 \pi
\end{equation}
are best suited for dealing with mappings of planar annuli. 
For a general Sobolev mapping $h$ we have the formula
\begin{equation}\label{polarjaco}
J_h(z)=J(z,h)= \frac{\textnormal{Im}\, (\overline{h_s} \, {h_\theta} )}{s} \leqslant \frac{|h_\theta |  \, |h_s|}{s} \, , \qquad s=\abs{z}\, . 
\end{equation}
We shall make use of the following free Lagrangians.
\begin{enumerate}[(i)]
\item{Pullback of a  form   in $\mathbb A^\ast$ via a given mapping $h \in \Hp{2}$;
\begin{equation*}
L(z,h,Dh)  = N(|h|)\, J(z,h) \, , \hskip0.5cm \textnormal{ where } \; N \in {\mathscr L}^1 (r_\ast , R_\ast)
\end{equation*} 
Thus, for all $h \in \Hp{2} $ we have
\begin{equation}\label{FL1}
\int_{\mathbb A} L(z,h,Dh)\, \dtext z = \int_{\mathbb A^\ast} N(|y|)\, \dtext y = 2 \pi \int_{r_\ast}^{R_\ast} N(G) \, G \, \dtext G
\end{equation}}
\item{A radial free Lagrangian
\begin{equation*}
L(z,h,Dh)\, \dtext z = A \big(|h|\big) \frac{|h|_s}{|z|}\, \dtext z\, , \hskip0.5cm \textnormal{ where } \; A \in \mathscr L^1 (r_\ast , R_\ast)
\end{equation*}
Thus, for all $h \in \Hp{2} $ we have
 \begin{equation}\label{FL2}
\int_{\mathbb A} L(z,h,Dh)\, \dtext z = 2 \pi \int_r^R A \big(|h|\big) \frac{\partial |h|}{\partial s}\, \dtext s = 2 \pi \int_{r_\ast}^{R_\ast} A(\tau) \, \dtext \tau
\end{equation}
}
\item{A tangential free Lagrangian
\begin{equation*}
L(z,h,Dh) = B \big( |z|\big) \textnormal{Im} \frac{h_\theta}{h} \, , \hskip0.5cm \textnormal{ where } \; B \in {\mathscr L}^1 (r , R)
\end{equation*}
Thus, for all $h \in \Hp{2} $ we have
\begin{equation}\label{FL4}
\int_{\mathbb \A} L(z,h,Dh)\, \dtext z = \int_r^R B(t) \left(\int_{|x|=t} \frac{\partial \textnormal{Arg}\, h}{\partial \theta} \, \dtext \theta\right) \dtext t =   2 \pi \int_{r}^{R} B(t) \, \dtext t
\end{equation}}
\item{Let $C(s, G)$, $r \le s \le R$, $r_\ast \le G \le R_\ast$ be a nonnegative $\mathscr C^1$-smooth function. The following differential expression is a free Lagrangian
\begin{equation*}
L(z,h,Dh) = (2C+G C_G)J(z,h)+  C_s  \frac{\abs{h}^2}{s}\im \frac{h_\theta}{h}\, , \quad s=\abs{z} \textnormal{ and } G=\abs{h} \, . 
\end{equation*}
For $h \in \Hp{2} $ we have
\begin{equation}\label{FL5}
\int_{\A} L(z,h,Dh)\, \dtext z  = \frac{\pi}{2} \left[ R_\ast^2 C(R,R_\ast)-r_\ast^2 C(r,r_\ast)  \right]
\end{equation}
}
\end{enumerate}
All of these Free-Lagrangians were introduces by Iwaniec and Onninen. The ones in (i)-(iii)  appeared first time in~\cite{IOan} and the last one (iv) in their forthcoming book~\cite{IOb}.

Let $\Omega \subset \mathbb C$ be a bounded Jordan domain with rectifiable boundary $\partial \Omega$. Then the familiar geometric form of the isoperimetric inequality reads as
\[\abs{\Omega} \le \frac{1}{4 \pi} [\ell (\partial \Omega)]^2\]
where $\abs{\Omega}$ is the area of $\Omega$ and $\ell (\partial \Omega)$ is the length of $\partial \Omega$. We denote the ball centered at the origin with radius $r$ by $B_r$. First, suppose that $f \colon B_R \to \C$ lies in the Sobolev class $W^{1,2} (B_R, \C)$. Partly using the polar coordinates $z=se^{i\theta}$, $s<R$ and $\theta \in [0, 2 \pi)$ we formulate the integral form of the isoperimetric inequality, see \cite{MUL,RESH}: 
\begin{equation}\label{eq:intiso}
\left| \int_{B_s} J(z, f) \, \dtext z \right| \le \frac{1}{4\pi} \left( \int_0^{2\pi} \abs{f_\theta}  \, \dtext \theta\right)^2 \quad \textnormal{ for almost every } s<R\, . 
\end{equation}
Second, suppose   that $h \colon \A \onto \A^\ast$  is  an   orientation-preserving diffeomorphism  which also preserves the order of boundary components. Then applying~\eqref{eq:intiso}  we obtain 
\begin{equation}\label{eq:isoapply}
\int_0^{2\pi} \im (\overline{h} h_\theta )\, \dtext \theta \le \frac{1}{2 \pi } \left( \int_0^{2\pi}  \abs{h_\theta} \, \dtext \theta \right)^2 \quad \textnormal{ for } r<s<R \, . 
\end{equation}
Indeed, fix $s\in (r,R)$. To simplify the notation we assume, without loss of generality, that $s=1$. Considering the differomorphism $h \colon \A \onto \A^\ast$ restricted to the unit sphere $S_1= \partial B_1$ and then extending this restricted mapping to the ball $B_R$ in the radial manner, namely,
\[f \colon B_R \to \C \, , \qquad f(z)= \abs{z}^2 h(z/\abs{z})\, . \]
Note $h(x)=f(x)$ on $S_1$ and $f$ is continuous differentiable on $B_R$. The isoperimetric inequality~\eqref{eq:intiso} yields,
\[ \int_{B_1} J(z, f) \, \dtext z \le \frac{1}{4\pi} \left( \int_0^{2\pi} \abs{f_\theta}  \, \dtext \theta\right)^2 = \frac{1}{4\pi} \left( \int_0^{2\pi} \abs{h_\theta}  \, \dtext \theta\right)^2\]
The Jacobian is the most known example of null Lagrangians and by Green's theorem  we have
\[ \int_{B_1} J(z, f) \, \dtext z  = \frac{1}{2} \int_0^{2\pi} \im (\overline{f} f_\theta )\, \dtext \theta =   \frac{1}{2} \int_0^{2\pi} \im ( \overline{h} h_\theta) \, \dtext \theta \, . \]
Therefore, the claimed version of the isoperimetric inequality~\eqref{eq:isoapply} follows.

\section{Proof of Theorem \ref{increasingThm}}
\begin{proof} \textbf{Case 1.} Assume that $R_*/r_* \geq m_\lambda(R/r)$.\\\\
The diffeomorphisms from $\A$ onto $\A^\ast$ are dense in $\Hp{2}$~\cite{IKOap}. Therefore, we can equivalently replace the admissible homeomorphism in~\eqref{eq:mainq} by diffeomorphisms. Precisely,
\[ \underset{\Hp{2}}{\inf} \mathbb E_\lambda [h] = \underset{\Diff (\A, \A^\ast)}{\inf} \mathbb E_\lambda [h]  \]
where  $\Diff (\A, \A^\ast)$ is for the class of orientation preserving diffeomorphisms from $\A$ onto $\A^\ast$ which also preserve the order of the boundary components.

Let $h\in  \Diff (\A, \A^\ast)$. We write $s=\abs{z}\in [r,R]$ and the weighted Dirichlet energy of $h$ in polar coordinates as follows
\begin{equation}\label{eq:blahblah}
\ee_\lambda[h] = \int_0^{2\pi} \int_r^R \frac{\lambda}{s} |h_\theta|^2 + s \lambda |h_s|^2 \, ds \, d\theta .
\end{equation}
Let us also denote by $W(Dh) = \frac{\lambda}{s} |h_\theta|^2 + s \lambda |h_s|^2$ the expression under the integral, which we will now estimate from below. Recall that $h_0$ denotes the minimizer among the radial mappings from $\A$ onto $\A^\ast$, which is a homeomorphism in this case. At this point we already remark that in all of the forthcoming estimates equality will hold for $h = h_0$, and $h_0$ will also be the only map for which there is equality in every estimate. Let us start by defining the expression
\[\tau = \Phi - \frac{c}{H},\]
where $c$ is a constant to be determined. Then our claim is that\\\\
\textbf{Claim 1.} If $R_*/r_* \geq m_\lambda(R/r)$, then the constant $c \geq 0$ may be chosen so that both $\tau \geq 0$ and $\dot{\tau} \geq 0$. Furthermore, the expressions $\lambda/s - \dot{\tau}$ and $s\lambda - c/\dot{H}$ are also nonnegative and we have the identity
\begin{equation}\label{tauIdentity}\left(\frac{\lambda}{s} - \dot{\tau}\right)\left(s\lambda - \frac{c}{\dot{H}}\right) = \tau^2.\end{equation}
After this claim is proven, our estimates for the expression $W$ proceed as follows
\begin{align*}
W &= \dot{\tau} |h_\theta|^2 + \left(\frac{\lambda}{s} - \dot{\tau}\right) |h_\theta|^2 + \left(s\lambda - \frac{c}{\dot{H}}\right)|h_s|^2 + c \frac{|h_s|^2}{\dot{H}}
\\&\geq \dot{\tau} |h_\theta|^2 + 2 \sqrt{\left(\frac{\lambda}{s} - \dot{\tau}\right)\left(s\lambda - \frac{c}{\dot{H}}\right)}|h_\theta||h_s| + c \frac{|h_s|^2}{\dot{H}}
\\&= \dot{\tau} |h_\theta|^2 + 2 \tau |h_\theta||h_s| + c \frac{|h_s|^2}{\dot{H}}.
\end{align*}
Here we applied to the elementary inequality $2ab \le a^2+b^2$ for real numbers $a,b$ and the identity~\eqref{tauIdentity}.
By a simple application of Cauchy-Schwartz, we obtain that
\begin{equation}\label{Westim1}\int_r^R \frac{|h_s|^2}{\dot{H}} ds \geq \frac{\left(\int_r^R h_s ds \right)^2}{\int_r^R \dot{H} ds} = R^* - r^*.\end{equation}
Here the equality is attained exactly for $h = h_0$, and we see that the right hand side does not depend on the choice of the map $h$. Next, we apply the isoperimetric inequality~\eqref{eq:isoapply} and the Cauchy-Schwartz inequality to find that
\begin{equation}\label{Westim2}\dot{\tau} \int_0^{2\pi} |h_\theta|^2 d\theta \geq \dot{\tau} \int_0^{2\pi}\im \left(\bar{h}h_\theta\right) d\theta.\end{equation}
Applying \eqref{Westim2}, we find that
\begin{equation}\label{Westim3}
\int_0^{2\pi} \int_r^R \dot{\tau} |h_\theta|^2 + 2 \tau |h_\theta||h_s| ds d\theta \geq \int_0^{2\pi} \int_r^R \dot{\tau} \im \left(\bar{h}h_\theta\right) + 2 \tau |h_\theta||h_s| ds d\theta.
\end{equation}
We are now in a position to apply the free Lagrangian (iv) from Section \ref{FLsection} with $C(s,G) = \tau(s)$ to the right hand side of~\eqref{Westim3}. The related equation \eqref{FL5} allows us to find the bound
\begin{align}\label{Westim4}
\int_0^{2\pi} & \int_r^R \dot{\tau} \im \left(\bar{h}h_\theta\right) + 2 \tau |h_\theta||h_s| ds d\theta \geq
\int_{\A} \dot{\tau}  \im \left(\bar{h}\frac{h_\theta}{s}\right) + 2\tau J_h\, dz
\\\nonumber &= \tau(R) 2 \pi R_*^2 - \tau(r) 2 \pi r_*^2,
\end{align}
which is independent of $h$. Combining the estimates \eqref{Westim1}-\eqref{Westim4}, we find the required lower bound for $\ee_\lambda[h]$.
\begin{align*}\ee_{\lambda}[h] &\geq \tau(R) 2\pi R_*^2 - \tau(r) 2\pi r_*^2 + 2\pi c (R_* - r_*)
\\&= 2\pi R_*^2 \Phi(R) - 2\pi r_*^2 \Phi(r) + 2\pi c \left(\frac{r_*^2}{H(r_*)} - \frac{R_*^2}{H(R_*)}\right) +  2\pi c (R_* - r_*)
\\&= 2\pi R_*^2 \Phi(R) - 2\pi r_*^2 \Phi(r)
\\&= \ee_{\lambda}[h_0]. 
\end{align*}
Let us now prove Claim 1. We start by verifying the identity \eqref{tauIdentity}. Here we make use of the equations \eqref{hdiff} and \eqref{phidiff}.
\begin{align*}
\left(\frac{\lambda}{s} - \dot{\tau}\right)\left(s\lambda - \frac{c}{\dot{H}}\right) 
&= \left(\frac{\lambda}{s} - \dot{\Phi} - \frac{c \dot{H}}{H^2}\right)\frac{H}{\dot{H}}\left(s\lambda \frac{\dot{H}}{H} - \frac{c}{H}\right)
\\&= \left(\frac{H}{\dot{H}}\left(\frac{\lambda}{s} - \dot{\Phi}\right) - \frac{c}{H}\right)\left(\Phi - \frac{c}{H}\right)
\\&= \left(\frac{s\lambda}{\Phi}\left(\frac{\lambda}{s} - \dot{\Phi}\right) - \frac{c}{H}\right)\tau
\\&= \left(\frac{1}{\Phi}\left(\lambda^2 - s\lambda \dot{\Phi}\right) - \frac{c}{H}\right)\tau
\\&= \left(\frac{1}{\Phi} \Phi^2 - \frac{c}{H}\right)\tau
\\&= \tau^2.
\end{align*}
Let us now choose the constant $c \geq 0$. We must choose this constant in such a way that both of the inequalities
\[\Phi -\frac{c}{H} \geq 0 \qquad \text{ and } \qquad \dot{\Phi} + \frac{c \dot{H}}{H^2} \geq 0\]
hold on every point of the interval $[r,R]$. Since $\dot{\Phi} = (\lambda^2 - \Phi^2)/(s\lambda)$ and $\dot{H}/H = \Phi/(s\lambda)$, we may transform these two inequalities into the following
\begin{equation}\label{cInequality}
H \Phi \geq c \geq \frac{H}{\Phi}\left(\Phi^2 - \lambda^2\right).
\end{equation}
Let us now make a couple of observations.\\\\
\emph{Observation 2.} The function $s \mapsto H(s) \Phi(s)$ is nondecreasing.\\\\
\emph{Proof.} By computation,
\begin{equation}\label{Obs1computation}
\frac{d}{ds} \left(H\Phi\right) = \dot{H}\Phi + H \dot{\Phi} = H \left(\frac{\Phi^2}{s\lambda} + \frac{\lambda^2 - \Phi^2}{s\lambda}\right) = \frac{\lambda H}{s} \geq 0.
\end{equation}
This observation shows that to satisfy the inequalities in \eqref{cInequality} we may as well choose $c = H(r) \Phi(r)$, as then the first inequality is always satisfied.\\\\
\emph{Observation 3.} Suppose that $\lambda$ is increasing. If at some point $s_0$ it holds that $\Phi(s_0) \leq \lambda(s_0)$, then $\Phi(s) \leq \lambda(s)$ for every $s \geq s_0$.\\\\
\emph{Proof.} The first sentence follows directly from \eqref{phidiff}. For the second part, observe that if it would hold that $\Phi(s_1) > \lambda(s_1)$ for some $s_1 > s_0$, then by continuity there would exist a point $s_2 \in [s_0,s_1)$ so that $\Phi(s_2) = \lambda(s_2)$ and $\Phi(s) \geq \lambda(s)$ for all $s \in [s_2,s_1]$. Since $\lambda$ is nondecreasing, we have $\Phi(s_1) > \lambda(s_1) \geq \lambda(s_2) = \Phi(s_2)$. By the mean value theorem we must have that $\dot{\Phi}(s_3) > 0$ for some point $s_3 \in [s_2,s_1]$, but this is a contradiction with Observation 1 since $\Phi(s) \geq \lambda(s)$ on this interval. \\\\
Since $\Phi$ is defined as a solution of the ODE \eqref{phidiff}, this observation shows that there are only the two following possible scenarios. Either $\Phi(r) \leq \lambda(r)$, in which case $\Phi(s) \leq \lambda(s)$ everywhere. In this case the right hand side in \eqref{cInequality} is nonpositive so any constant $c \in [0,\Phi(r)H(r)]$ will do.

In the second case, $\Phi(r) > \lambda(r)$. In this case $\Phi$ starts out as decreasing, and may hit $\lambda$ at some point $s_0 \in (r,R]$. If such a $s_0$ exists then $\Phi(s) \leq \lambda(s)$ for all $s \in [s_0,R]$ and the right inequality in \eqref{cInequality} holds on this part of the interval $[r,R]$. In any case, it is enough to show that the right inequality in \eqref{cInequality} holds on an interval $[r,s_0)$ on which $\Phi$ is decreasing. But this is an easy consequence of the following.\\\\
\emph{Observation 4.} At the points where $\Phi$ is decreasing, the expression $\frac{H}{\Phi}\left(\Phi^2 - \lambda^2\right)$ is also decreasing.\\\\
\emph{Proof.} We make a direct computation. Here we also reuse the computation \eqref{Obs1computation}.
\begin{align*}
\frac{d}{ds}\left(\frac{H}{\Phi}\left(\Phi^2 - \lambda^2\right)\right) &= \frac{\lambda H}{s} - \frac{d}{ds} \frac{H \lambda^2}{\Phi}
\\&= \frac{\lambda H}{s} - \frac{\dot{H} \lambda^2}{\Phi} - \frac{H 2\dot{\lambda}\lambda}{\Phi} + \frac{\dot{\Phi} H \lambda^2}{\Phi^2}
\\&= - \frac{H 2\dot{\lambda}\lambda}{\Phi} + \frac{\dot{\Phi} H \lambda^2}{\Phi^2}
\\&= \frac{\lambda^4 H}{\Phi^2}\left(- \frac{2 \Phi \dot{\lambda}}{\lambda^3} + \frac{\dot{\Phi}}{\lambda^2}\right)
\\&= \frac{\lambda^4 H}{\Phi^2}\, \frac{d}{ds} \left(\frac{\Phi}{\lambda^2}\right).
\end{align*}
Since $\Phi$ is decreasing and $\lambda$ is nondecreasing, the expression $\Phi/\lambda^2$ is decreasing. Thus the last expression above is negative, and we have proved the claim.\\\\
Via Observation 4, we now find that for $s$ on the interval $[r,s_0)$ we have
\[\frac{H(s)}{\Phi(s)}\left(\Phi^2(s) - \lambda^2(s)\right) < \frac{H(r)}{\Phi(r)}\left(\Phi^2(r) - \lambda^2(r)\right) < H(r)\Phi(r) = c.\]
This proves the inequality \eqref{Obs1computation} for the choice of $c = H(r)\Phi(r)$.\\\\
Returning to the statement of Claim 1, we must still verify the nonnegativity of the expressions $\lambda/s - \dot{\tau}$ and $s\lambda - c/\dot{H}$. But this easily follows from the nonnegativity of $\tau, H$ and $\dot{H}$ as well as the identity
\[\tau = \frac{\dot{H}}{H} \left(s\lambda - \frac{c}{\dot{H}}\right) = \frac{\dot{H}}{H} \frac{\tau^2}{\lambda/s - \dot{\tau}},\]
which was essentially verified in the proof of \eqref{tauIdentity}.\\\\
\textbf{Case 2.} Assume that $R_*/r_* < m_\lambda(R/r)$.\\\\
Recall that in this case there exists a radius $r_0 \in (r,R)$ so that $H(s) = r_*$ for all $s \in [r,r_0]$. For the corresponding radial minimizer $h_0\colon \A \to \overline{\A^\ast}$ the
part of the domain annulus near its inner boundary  collapses into the inner boundary of $\A^\ast$. Moreover, the function $\Phi$ is identically zero on $[r,r_0]$ and solves the equation \eqref{phidiff} only on $(r_0,R]$. This suggests that we should estimate  the integral~\eqref{eq:blahblah} into two separate parts.

On the interval $ [r,r_0]$, we apply the estimates
\begin{align*}
\int_0^{2\pi} \int_r^{r_0} \frac{\lambda}{s} |h_\theta|^2 + s \lambda |h_s|^2 \, ds \, d\theta 
&\geq \int_r^{r_0}  \frac{\lambda}{s} \int_0^{2\pi} |h_\theta(s e^{i\theta})|^2 \, d\theta \, ds
\\&\geq \int_r^{r_0}  \frac{\lambda}{s} \frac{1}{2\pi}\left(\int_0^{2\pi} |h_\theta(s e^{i\theta})| \, d\theta\right)^2 \, ds
\\&\geq \int_r^{r_0}  \frac{\lambda}{s} \frac{1}{2\pi}(2\pi r_*)^2 \, ds
\\&= 2\pi r_*^2 \int_r^{r_0}  \frac{\lambda(s)}{s}  ds,
\end{align*}
where the last inequality is due to the fact that the length of the image curve of $s e^{i\theta}$ with $\theta \in [0,2\pi)$ under $h$ is at least $2\pi r_*$. In particular, equality here holds exactly for $h = h_0$ since $h_0$ sends the annulus $\A(r,r_0)$ to the circle of radius $r_*$.

On the interval $(r_0,R]$ we apply the same estimates as in Case 1. However, in this case we may simply choose the constant $c$ appearing in Case 1 to be zero, as the fact that $\Phi(r_0) = 0$ implies that $\Phi(s) \leq \lambda(s)$ everywhere. Hence we  have that $\tau = \Phi$. This results in the estimate
\begin{align*}
\int_0^{2\pi}& \int_{r_0}^R \frac{\lambda}{s} |h_\theta|^2 + s \lambda |h_s|^2 \, ds \, d\theta
\\&\geq  \int_0^{2\pi} \int_{r_0}^R \dot{\Phi} \im \left(\bar{h} h_\theta \right) + 2 \Phi |h_\theta||h_s| ds d\theta
\\& \geq \int_0^{2\pi}  \Phi(R) \im\left[ \bar{h}(R e^{i\theta}) h_\theta(R e^{i\theta}) \right] - \Phi(r_0) \im\left[ \bar{h}(r_0 e^{i\theta}) h_\theta(r_0 e^{i\theta}) \right] d\theta
\\&= \Phi(R) \int_0^{2\pi} \im\left[ \bar{h}(R e^{i\theta}) h_\theta(R e^{i\theta}) \right] d\theta
\\&= 2 \pi R_*^2 \Phi(R)
\end{align*}
Combined, we have that
\[\ee_{\lambda}[h] \geq 2 \pi R_*^2 \Phi(R) + 2\pi r_*^2 \int_r^{r_0}  \frac{\lambda(s)}{s}  ds = \ee_{\lambda}[h_0].\]
\end{proof}

\section{Proof of Theorem \ref{thinThm}}
\begin{proof} The proof of this theorem follows the same lines of arguments as the proof of Theorem \ref{increasingThm}. The only parts in the proof of Theorem \ref{increasingThm} where the fact that $\lambda$ is nondecreasing was used were
\begin{enumerate}
\item{To guarantee that if $\Phi(s_0) \leq \lambda(s_0)$ then $\Phi(s) \leq \lambda(s)$ for all $s > s_0$.}
\item{To deal with the estimates for the constant $c$ at the points $s$ where $\Phi(s) \geq \lambda(s)$.}
\end{enumerate}
Hence if we are somehow able to guarantee that $\Phi(s) \leq \lambda(s)$ for every point $s \in [r,R]$, then the proof of Theorem \ref{increasingThm} adapts to any positive continuous weight. However, if we recall the discussion in Section \ref{radialsection}, this is always possible to do by choosing a small enough initial value $\varphi_0$. Letting $\Phi_2$ denote the largest map $\Phi$ for which the inequality $\Phi(s) \leq \lambda(s)$ holds for every $s \in [r,R]$, we may define the function $g_\lambda$ by the formula
\begin{equation}\label{smallglambda}
\log g_\lambda(R/r) = \int_r^R \frac{\Phi_2(s)}{s \lambda(s)} \, ds.
\end{equation}
This definition guarantees that for any target $\A^*$ with $R_*/r_* \leq g_\lambda(R/r)$ the associated function $\Phi$ satisfies $\Phi(s) \leq \Phi_2(s) \leq \lambda(s)$ for every $s \in [r,R]$, which proves the fact that $\ee_\lambda[h]$ has a radial minimizer by the discussion above.

The fact that $\lambda$ is continuous and positive on $[r,R]$ implies that $\lambda$ is bounded from below by a positive constant, which in turn guarantees that the map $\Phi_2$ defined above is not identically zero. This also shows that $g_\lambda(x) > 1$ for every $x > 1$.
\end{proof}

\section{Proof of Theorem \ref{fixedThm}}
The proof of this theorem is based on the proof of Theorem 1.3 in \cite{KOpharm}. We begin by applying H\"older's inequality in the form
\[\int_\A |Dh|^2 \lambda dz \geq \frac{\left( \int_\A |Dh_0| |Dh|  \lambda  dz\right)^2}{\int_\A |Dh_0|^2  \lambda  dz},\]
where $h_0$ denotes the radial minimizer as defined in Section 2. Since equality holds here for $h = h_0$, it will be sufficient to estimate the quantity $\int_\A |Dh_0| |Dh|  \lambda  dz$ on the right hand side. Let $g(s)$ be a function to be determined, $0 \leq g(s) \leq 1$. Writing $|Dh|$ in polar coordinates and applying an elementary inequality gives
\[|Dh| = \sqrt{|h_s|^2 + \left|\frac{h_\theta}{s}\right|^2} \geq \sqrt{g(s)} |h_s| + \sqrt{1-g(s)}\frac{|h_\theta|}{s}.\]
We wish to find $g$ such that equality holds in this estimate for $h = h_0$. A short calculation gives
\[\frac{\sqrt{g(s)}}{\sqrt{1-g(s)}} = \frac{s \dot{H}(s)}{H(s)} = \frac{\Phi(s)}{\lambda(s)} \quad \Leftrightarrow \quad g(s) = \frac{\Phi^2(s)}{\Phi^2(s) + \lambda^2(s)}.\]
We hence obtain the estimate
\[\int_\A |Dh_0| |Dh|  \lambda  dz \geq \int_0^{2\pi}\int_r^R \rho_1(s) |h_s| + \rho_2(s) |h_\theta| ds d\theta\]
where one may compute that the coefficients simplify to
\[\rho_1(s) = \Phi(s)H(s) \qquad \text{ and } \qquad \rho_2(s) = \frac{ \lambda(s) H(s)}{s}.\]
From the computation \eqref{Obs1computation} we may see that these coefficients satisfy the equality $\dot{\rho_1}(s) = \rho_2(s)$. The rest of the proof is exactly the same as in \cite{KOpharm}, following from the part of the proof of Theorem 1.3 after a similar equality was established. The key assumption, i.e. that the mapping $h$ is fixed on the outer boundary, is utilized in this part of the proof just before Lemma 3.1.

\end{document}